\documentclass{amsart}
\usepackage{amssymb}

\newtheorem{thm}{Theorem}[section]
\newtheorem{cor}[thm]{Corollary}
\theoremstyle{definition}
\newtheorem{dfn}[thm]{Definition}
\newtheorem{ex}[thm]{Example}
\theoremstyle{remark}
\newtheorem{rmk}[thm]{Remark}
\numberwithin{equation}{section}

\begin{document}

\title{Level Algebras with Bad Properties}

\author{Mats Boij}
\address{Department of Mathematics, Royal Institute of Technology, S-100 44 Stockholm, Sweden}
\email{boij@math.kth.se}

\author{Fabrizio Zanello}
\address{Department of Mathematics, Royal Institute of Technology, S-100 44 Stockholm, Sweden}
\email{zanello@math.kth.se}
\thanks{The second author is funded by the G\"oran Gustafsson Foundation}

\date{}

\begin{abstract}
This paper can be seen as a continuation of the works contained in the recent preprints \cite{Za}, of the second author, and \cite{Mi}, of Juan Migliore. Our results are:

1). There exist codimension three artinian level algebras of type two which do not enjoy the Weak
 Lefschetz Property (WLP). In fact, for $e\gg 0$, we will construct a codimension three, type two $h$-vector of socle degree $e$ such that {\em all}
the level algebras with that $h$-vector do not have the WLP. We will also describe the family of those algebras and compute its dimension, for each $e\gg 
0$.

2). There exist reduced level sets of points in ${\mathbf P}^3$ of type two whose artinian reductions all fail to have the WLP. Indeed, the examples constructed here have the same $h$-vectors we mentioned in 1). 

3). For any integer $r\geq 3$, there exist non-unimodal monomial artinian level algebras of codimension $r$. As an immediate consequence of this result, we obtain another proof of the fact (first shown by Migliore in \cite[Theorem 4.3]{Mi}) that, for any $r\geq 3$, there exist reduced level sets of points in ${\mathbf P}^r$ whose artinian reductions are non-unimodal.
\end{abstract}
\maketitle
\section{Introduction}
The study of the properties of level algebras, both in the artinian and in the higher dimensional case, is an intriguing issue in commutative algebra, not only because of the intrinsic importance of these algebras, but also due to the relevance of their applications to several other areas of mathematics, such as algebraic geometry, invariant theory, combinatorics, etc. 
There is an excellent broad overview of level algebras in the memoir \cite{GHMS}, and we refer to it for the most comprehensive bibliography up to the year 2003. However, although a remarkable amount of research on level algebras has been performed over the last few years, several interesting problems are still open. The purpose of this note is to provide the solution to three of those problems.

We will first focus on the Weak Lefschetz Property (WLP, in brief; see below for the definition), which is a fundamental and very natural property of artinian algebras, and has recently received much attention (e.g., see \cite{HMNW} and \cite{MM} and their bibliographies). In codimension $r\geq 4$, it is known that the WLP may fail even for Gorenstein algebras (see Ikeda's \cite[Example 4.4]{Ik}, and the first author's \cite[Theorem 3.6]{Bo2}, for $r=4$; see also \cite[Example 4.3]{St},  \cite{BI}, \cite{BL} and \cite{Bo1}, which supply examples of Gorenstein algebras of codimension $r\geq 5$ that are non-unimodal, and therefore, {\it a fortiori}, without the WLP).

The second author (see \cite[Corollary 4 and Proposition 8]{Za}) has recently shown that there even exist codimension 3 artinian level algebras which do not enjoy the WLP, but could only provide examples having type greater than or equal to 3. (The idea, at least for type three, was to use an inverse system generated by three forms of the same degree, two suitably chosen in two variables, $y_1$, $y_2$, and the third being generic in variables $y_2$, $y_3$.) Therefore, only the problem as to whether there exist Gorenstein or type 2 level algebras of codimension 3 without the WLP remained open. As for the Gorenstein case, so far the most important result is \cite[Corollary 2.4]{HMNW},  which proves, with a beautiful geometric argument, that all complete intersections of codimension 3 enjoy the WLP.

The first result of this paper is the solution to the type 2 case: indeed, we will show that, for any $e\gg 0$, there exist codimension 3, type 2 level algebras of socle degree $e$ which do not enjoy the WLP. Our examples are also surprisingly simple. In fact, the inverse system modules we will construct are generated by one monomial and one binomial. Actually, we will prove more: we will determine the $h$-vectors $h^{(e)}$ of the level algebras without the WLP mentioned above and show that, in fact, {\it all} the level algebras having those $h$-vectors do not enjoy the WLP. For each $h^{(e)}$, we will also completely describe the family of the level algebras with $h$-vector $h^{(e)}$, and show that it is an irreducible variety of dimension $3(e+1)$.

The recent preprint of Migliore's, \cite{Mi}, has extended some of the results of \cite{Za} to (one-dimensional) level algebras of type $t\geq 3$ of reduced sets of points of the projective space ${\mathbf P}^3$. The second result of this note, obtained by employing a technique rather similar to that of Migliore, is the construction of type 2, reduced level sets of points of ${\mathbf P}^3$ whose artinian reductions all fail to have the WLP. Indeed, these examples have the same $h$-vectors as those supplied in our first result mentioned above. Thus, as far as the WLP is concerned, both for artinian and for one-dimensional codimension 3 level algebras the non-existence problem remains now open only in the Gorenstein case.

In \cite[Theorem 3]{Za}, the second author has also proved the existence of codimension 3 artinian level algebras having a non-unimodal $h$-vector. (These examples all have a relatively large type, and for small types the unimodality problem is still open.) However, none of his algebras were monomial, and the problem as to whether there exist non-unimodal monomial level algebras (of any codimension) remained unsolved.

The third main result of this note is the solution to that problem: we will show that, for every $r\geq 3$, there exist codimension $r$ monomial artinian level algebras having a non-unimodal $h$-vector. As an immediate corollary, since monomial artinian ideals lift to ideals of points, we have another proof of the existence of reduced level sets of points in ${\mathbf P}^r$ (for every $r\geq 3$) whose artinian reductions are non-unimodal (this result was first shown in \cite[Theorem 4.3]{Mi}).

\section{Setup}
Let us now introduce the main definitions we will need in this note. We consider standard graded algebras $A=R/I$, where $R=k[x_1,...,x_r]$, $I$ is a homogeneous ideal of $R$, $k$ is a field of characteristic zero and the $x_i$'s all have degree 1. Our algebras will be either artinian or Cohen-Macaulay (CM) (mostly of dimension $d=1$, like in the case of points in the projective space).

When $A$ is artinian, the {\it $h$-vector} of $A$ is $h(A)=h=(h_0,h_1,...,h_e)$, where $h_i=\dim_k A_i$ and $e$ is the last index such that $\dim_k A_e>0$. Since we may suppose, without loss of generality, that $I$ does not contain non-zero forms of degree 1, the codimension of $A$ is $r=h_1$. In general, when $A$ is CM of dimension $d$, its $h$-vector is that of any of its artinian reductions (i.e., artinian quotients of $A$ by a regular sequence), and its codimension is $r-d$.

In the artinian case, the {\it socle} of $A$ is the annihilator of the maximal homogeneous ideal $\overline{m}=(\overline{x_1},...,\overline{x_r})\subseteq A$, namely $soc(A)=\lbrace a\in A {\ } \mid {\ } a\overline{m}=0\rbrace $. Since $soc(A)$ is a homogeneous ideal, we define the {\it socle-vector} of $A$ as $s(A)=s=(s_0,s_1,...,s_e)$, where $s_i=\dim_k soc(A)_i$. Note that $h_0=1$, $s_0=0$ and $s_e=h_e>0$. The integer $e$ is called the {\it socle degree} of $A$ (or of $h$). The {\it type} of the socle-vector $s$ (or of the algebra $A$) is type$(s)=\sum_{i=0}^es_i$.

If $s=(0,0,...,0,s_e=t)$, we say that the graded algebra $A$ is {\it level} (of type $t$). If, moreover, $t=1$, then $A$ is {\it Gorenstein}.

Since the graded Betti numbers of a minimal free resolution (MFR) of a CM algebra $A$ are preserved when we consider the MFR of any artinian reduction $A^{'}$ of $A$, and since the non-zero entries of the socle-vector of $A^{'}$ are given by the graded Betti numbers of the last module of its MFR, then we can define a CM algebra as {\it level of type $t$} when its artinian reductions are level of type $t$.

An artinian algebra $A=\oplus_{i=0}^e A_i$ is said to have the {\it Weak Lefschetz Property} ({\it WLP}) if there exists a linear form $l\in R$ such that, for all indices $i=0,1,...,e-1$, the multiplication map \lq \lq $\cdot l$" between the $k$-vector spaces $A_i$ and $A_{i+1}$ has maximal rank.

Let us now recall the main facts of the theory of {\it inverse systems}, or {\it Macaulay duality}, which will be a fundamental tool in this note. For a complete introduction, we refer the reader to \cite{Ge} and \cite{IK}.

Let $S=k[y_1,...,y_r]$, and consider $S$ as a graded $R$-module where the action of $x_i$ on $S$ is partial differentiation with respect to $y_i$.

There is a one-to-one correspondence between graded artinian algebras $R/I$ and finitely generated graded $R$-submodules $M$ of $S$, where $I=Ann(M)$ is the annihilator of $M$ in $R$ and, conversely, $M=I^{-1}$ is the $R$-submodule of $S$ which is annihilated by $I$ (cf. \cite[Remark 1]{Ge}), p. 17).

If $R/I$ has socle-vector $s$, then $M$ is minimally generated by $s_i$ elements of degree $i$, for $i=1,2,...,e$, and the $h$-vector of $R/I$ is given by the number of linearly independent partial derivatives obtained in each degree by differentiating the generators of $M$ (cf. \cite[Remark 2]{Ge}), p. 17).

In particular, artinian level algebras of type $t$ and socle degree $e$ correspond to $R$-submodules of $S$ minimally generated by $t$ elements of degree $e$.

It is easy to see that the inverse system module $M$ of an ideal $I$ is generated by monomials if and only if $I$ is a monomial ideal (i.e., it is also generated by monomials).

\section{Non-WLP for type two - The artinian case}

We now turn our attention to the first result of this note: as we have mentioned above, for each $e\gg 0$, we will construct a codimension 3 $h$-vector of socle degree $e$, such that all the type 2 level algebras having this $h$-vector do not enjoy the WLP.

Let $R=k[x_1,x_2,x_3]$, $S=k[y_1,y_2,y_3]$, and let \lq \lq $\circ $" indicate the operation of differentiation on inverse system modules.

\begin{dfn} For every integer $e$, define the family ${\mathcal F}_e$ by
\begin{equation*}{\mathcal F}_e=\lbrace R/Ann(M) {\ } \mid {\ } M=(F,G)\subset S, \mbox{ where $\deg F=\deg G=e$ and}\end{equation*}\begin{equation*} l_1\circ F=(l_1l_2)\circ G=0 \mbox{ for some, not necessarily distinct, linear forms $l_1,l_2\in R$} \rbrace ,\end{equation*}
and define the $h$ vector $h^{(e)}=(h_0,h_1,\dots,h_e)$ by 
\begin{equation*}
h_d= \min \lbrace 2d+1,e+2,3(e-d)+2\rbrace, 
\end{equation*}
for $d=0,1,\dots,e$.
\end{dfn}

\begin{thm}\label{thm1}
The generic element of $\mathcal F_e$ has $h$-vector $h^{(e)}$, and for $e\ge 6$ all level algebras with $h$-vector $h^{(e)}$ are in $\mathcal F_e$, which is an irreducible variety of dimension $3(e+1)$.

Moreover, for odd $e\ge 9$ and even $e\ge 12$, any level algebra with $h$-vector $h^{(e)}$ does not satisfy the Weak Lefschetz Property. 
\end{thm}

\begin{proof}
We start by finding an element in $\mathcal F_e$ with $h$-vector $h^{(e)}$. Namely, we want to show that 
\begin{equation*}
M=<y_1^{[2e/3]}y_3^{[e/3]},y_2^{\lceil e/2\rceil }y_3^{\lfloor e/2\rfloor }+y_1^{[e/3]}y_3^{[2e/3]}>\subset S
\end{equation*} 
is the inverse system module associated to a type 2 level algebra $A$, which is clearly in $\mathcal F_e$, whose $h$-vector is $h^{(e)}$. Recall that $\lceil a\rceil $ indicates the smallest integer greater than or equal to $a$, $\lfloor a\rfloor $ the largest integer less that or equal to $a$, and $[a]$ the integer nearest to $a$, assuming that this can be uniquely determined. There are six different cases to consider, one for each equivalence class of $e$ modulo 6. We will only prove the result for $e \equiv 4$ (mod 6); the arguments for the other cases are entirely similar and therefore will be omitted.

Hence let $e=6n+4$. Thus, 
\begin{equation*}
M=(F,G)=(y_1^{4n+3}y_3^{2n+1},y_2^{3n+2}y_3^{3n+2}+y_1^{2n+1}y_3^{4n+3}),
\end{equation*} 
and we want to show that the $h$-vector of $R/Ann(M)$ is $(1,3,5,...,6n+5,6n+6,6n+6,...,6n+6,6n+5,6n+2,...,11,8,5,2)$, where the first $6n+5$ is in degree $3n+2$, and $6n+6$ appears $n$ times (from degree $3n+3$ to degree $4n+2$).

It is easy to see that, for each $i=1,2,...,2n+1$, the $i$-th (linearly independent partial) derivatives of $F$ are $i+1$ (i.e., they are as many as they can possibly be for a form in two variables), and the derivatives of $G$ are $2i+1$ (since $y_1y_2$ does not appear in the summands of $G$). Furthermore, the first $2n+1$ derivatives of $F$ and $G$ are clearly disjoint. Hence, for $1\leq i\leq 2n+1$, $h_{e-i}=(i+1)+(2i+1)=3i+2$. This proves that $h$ ends with $(...,6n+5,6n+2,...,8,5,2)$ from degree $4n+3$ on.

It is immediate to see that $h_1=3$ and $h_2=5$ (since $x_1x_2$ is the only degree 2 monomial annihilating both $F$ and $G$). Hence, by Macaulay's theorem and the properties of the binomial expansion, it is easy to see that the equality $h_d=2d+1$ for $d\leq 3n+2$ follows once we show that $h_{3n+3}=6n+6$. Therefore, it remains to prove that $h_{3n+3}=h_{3n+4}=...=h_{4n+2}=6n+6$.

A simple computation shows that $Ann(M)_{4n+2}$ is generated (as a $k$-vector space) by monomials, with $A_{4n+2}$ spanned by the classes of the following two sets of monomials: \begin{equation*}T_1=\lbrace x_1^{4n+2},x_1^{4n+1}x_3, ..., x_3^{4n+2}\rbrace \end{equation*}
(i.e., all the $4n+3$ degree $4n+2$ monomials of $R$ in variables $x_1$ and $x_3$ only), and \begin{equation*}T_2=\lbrace x_2^{3n+2}x_3^n,x_2^{3n+1}x_3^{n+1},...,x_2^nx_3^{3n+2}\rbrace ,\end{equation*}
which has cardinality $3n+2-n+1=2n+3$. Thus, $h_{4n+2}=(4n+3)+(2n+3)=6n+6$, as desired.

Now, for each $i=1,2,...,n-1$, the monomials generating $A_{4n+2-i}$ are of course obtained by differentiating the monomials of $T_1$ and $T_2$. Those coming from $T_1$ are exactly the $4n+3-i$ monomials of degree $4n+2-i$ in variables $x_1$ and $x_3$ only, while those obtained from $T_3$ are $x_2^{3n+2}x_3^{n-i},x_2^{3n+1}x_3^{n-i+1},...,x_2^{n-i}x_3^{3n+2}$, which are $2n+3+i$. Hence, for $i=1,2,...,n-1$, $h_{4n+2-i}=(4n+3-i)+(2n+3+i)=6n+6$, as we wanted to show.

We now proceed to prove that $h^{(e)}$ is an upper bound for the $h$-vector of any member of the family $\mathcal F_e$. For this purpose, let $A=R/I$ be an element of $\mathcal F_e$, where $I=Ann(F,G)$ for some $F$ and $G$ in $S_e$.

Since $l_1 l_2\circ F=l_1l_2\circ G=0$, for some linear forms $l_1$ and $l_2$, we have that $l_1l_2\in I$. Thus the $h$-vector of $A$ is bounded from above by the Hilbert function of $R/(l_1l_2)$, i.e., $2d+1$.  

The inverse system module $I^{-1}$ is generated by $F$ and $G$ with the two relations $l_1\circ F=0$ and $l_1l_2\circ G=0$. Thus we get the inequality
\begin{equation*}
h(A)_{e-d} \le (d+1) + (2d+1) = 3d+2,\qquad \hbox{for $d=0,1,\dots,e$.}
\end{equation*}

We now have to establish the inequality $h(A)_d \le e+2$. By a change of variables we may assume that $l_1=x_1$. The derivatives of $F$ and $l_2\circ G$ with respect to $y_1$ are all zero. This means that the derivatives of degree  $e-d$ of $F$ and $l_2\circ G$ live in a space of dimension $e-d+1$. What remains are the derivatives of $G$ obtained by operating with forms of degree $d$ in $x_2$ and $x_3$, and the dimension of their span is bounded from above by $d+1$. Thus, in total, the $h$-vector of $A$ in degree $e-d$ is bounded from above by $(e-d+1) + (d+1)= e+2$. 

We now prove that any level algebra with $h$-vector $h^{(e)}$, for $e\ge 6$, is in the family $\mathcal F_e$. Let $A=R/I$ be any artinian level algebra with $h(A)=h^{(e)}$. By the value of the $h$-vector in degree two, we see that there is a quadratic form $q$ in the ideal $I$. Let $F$ and $G$ be the generators of the inverse system module $I^{-1}$. By the value of the $h$-vector in degree $e-1$, we have a linear relation among $F$ and $G$, i.e., there are linear forms $l_1$ and $l_2$ such that $l_1\circ F = l_2 \circ G$. Suppose that both of these forms are non-zero and linearly independent. Then we can find a form $H$ such that $F = l_2 \circ H$ and $G=l_1\circ H$. Since $e\ge 6$, we have that $h(A)_{e-2} = 8$ and we need eight linearly independent second derivatives of $F$ and $G$. There are at most nine possible such derivatives, since $F$ and $G$ are derivatives of $H$, and in addition we have $q\circ F = q\circ G =0$. Thus  $h(A)_{e-2}\le 7$, and we can conclude that the linear relation can be written as $l_1\circ F=0$. Hence we have that $x_1l_1\circ F = x_2l_1\circ F = x_3l_1\circ F=0$. Together with $q\circ F=q\circ G=0$, we get five quadratic relations among $F$ and $G$. Because of the $h$-vector, these relations must be linearly dependent, i.e., $q = l_1(ax_1+bx_2+cx_3) = l_1l_2$ for some linear form $l_2 = ax_1+bx_2+cx_3$. This proves that $A$ is in the family $\mathcal F_e$. 

Now, we can show that any element of $\mathcal F_e$ with Hilbert function $h^{(e)}$ does not satisfy the WLP for odd $e\ge 9$ and even $e\ge 12$.

Let $l$ be a generic linear form in $A$. Then the ideal $(l)$ is a cyclic submodule of $A$ and hence isomorphic to a quotient $B=A/(0:l)=R/(I:l)$. Since $(I:l)^{-1} = (l\circ F,l\circ G)$, we have that $B$ is an element of $\mathcal F_{e-1}$. Hence all entries of the $h$-vector of $B$ are bounded by $(e-1)+2=e+1$. This means that the Hilbert function of $(0:l)$ is at least $1$ in degrees where $h(A)_d=e+2$. Thus the multiplication by $l$ fails to be injective in these degrees. For odd $e\ge 9$ and even $e\ge 12$, we have that $h^{(e)}$ attains the value $e+2$ at least twice. Thus, in these cases, the algebra $A$ does not have the WLP.

We can prove that the family $\mathcal{F}_e$ is irreducible of dimension $3(e+1)$ in the following way. We map the family onto the space of reducible conics, $q=l_1l_2$, which has dimension four. The fibers are irreducible of dimension $(e +1) + (2e +1) -3 = 3e-1$, since the choice of $F$ and $G$ corresponds to $e+1$ and $2e+1$, respectively, but $(F,G)=(aF,bG+cF)$, whenever $a$ and $b$ are non-zero. Thus the family $\mathcal F_e$ is irreducible of dimension  $4+(3e-1) = 3(e+1)$, as we wanted to show.
\end{proof}

\section{Non-WLP for type two - The case of points}
We will now show that our first result can be ``lifted'' to homogeneous coordinate rings of points of the projective space $\mathbf P^3$ which are level of type two and whose artinian reductions all fail to satisfy the WLP. The idea is rather similar to that of Migliore (\cite{Mi}), which involved the use of sets of points on two different hyperplanes. In order to obtain the desired $h$-vector and type we need to be careful in the choice of the points on the hyperplanes. We can view the technique we will employ here as a basic double linkage, as in Migliore's work, although we will not use this description in the proof.

For any integer $e$, define the positive integers $c_0, c_1, c_2$ and $q_0,q_1,q_2,q_3$ such that: $c_0\leq c_1\leq c_2\leq c_0+1$, $c_0+c_1+c_2=e$, $|q_i-q_j|\leq 1$ for all $i$ and $j$, and
\begin{equation*}\left\{
\begin{array}{rl}
q_0+q_1 &=  c_0+1,\\
q_1+q_2 &=  c_1+1,\\
q_2+q_3 &=  c_2+1.\\
\end{array}
\right.\end{equation*}\indent
Equivalently, in a more compact fashion, we can write
\begin{equation*} c_i = \left\lfloor \frac{e+i}{3}\right\rfloor, \qquad \hbox{for $i=0,1,2$,} \end{equation*}
and
\begin{equation*} q_i =\left\lfloor \frac{e+i-1}{2}\right\rfloor - \left\lfloor  \frac{e+i-2}{3}\right\rfloor, \qquad \hbox{for $i=0,1,2,3$.}\end{equation*}\indent
Define the ideal $I\subset R=k[x_0,x_1,x_2,x_3]$ by 
\begin{equation*}\begin{array}{rl}
I &= (l_0l_1, l_0(C_0Q_2-C_1Q_0) ,l_0(C_1Q_3-C_2Q_1), \\
&C_0Q_1Q_2-C_1Q_0Q_1,C_1Q_2Q_3-C_2Q_1Q_2,C_0Q_2Q_3-C_2Q_0Q_1),
  \end{array}\end{equation*}
where $l_0,l_1$ are linear forms and $C_0,C_1,C_2,Q_0,Q_1,Q_2,Q_3$ are forms of degrees $c_0,c_1,c_2,q_0,q_1,q_2,q_3$.

\begin{thm}\label{thm2}
For a generic choice of the forms $l_0,l_1,C_0,C_1,C_2,Q_0,Q_1,Q_2$ and $Q_3$, the ideal $I\subset R$ defined above is the ideal of a reduced level set of points of ${\mathbf P}^3$ of type two, having $h$-vector $(1,h_1,...,h_e)$, where
\begin{equation*}h_d = \min \{2d+1,e+2,3(e-d)+2\}, \qquad \hbox{for $d=0,1,...,e$.}\end{equation*}\indent
By Theorem~\ref{thm1}, no artinian reduction of $A=R/I$ satisfies the WLP.
\end{thm}

\begin{proof}
Let us consider the ideal $I_1$ of $R$, generated by the maximal minors of the matrix
\begin{equation*}\left(
\begin{array}{ccc}
C_0&C_1&C_2\\
Q_0Q_1&Q_1Q_2&Q_2Q_3\\
\end{array}
\right),\end{equation*}
and the ideal $(f,g)$, where $f=C_1Q_3-C_2Q_1$ and $g=C_0Q_2-C_1Q_0$. The ideal $I_1$ is the ideal of a curve containing the complete intersection curve defined by $(f,g)$. The residual scheme is given by
\begin{equation*}I_1:(f,g)=(Q_0Q_1,Q_1Q_2,Q_2Q_3) = (Q_0,Q_2)\cap(Q_1,Q_2)\cap(Q_1,Q_3),\end{equation*}
that is the intersection of three complete intersection ideals.

For general forms $Q_0,Q_1,Q_2,Q_3$, each complete intersection ideal is radical, and so is their intersection $I_1:(f,g)$. The complete intersection ideal $(f,g)$ is also radical for a general choice of the above seven forms. To see this, we can choose a special set of forms and show that $(f,g)$ is radical for this choice. In order to do that, let $C_1=0$. Then $(f,g)=(C_0Q_2,C_2Q_1)=(C_0,C_2)\cap(C_0,Q_1)\cap(Q_2,C_2)\cap(Q_1,Q_2)$ is clearly radical if we take $C_1,C_2,Q_1,Q_2$ to be products of generic linear forms. Thus, since $I=((l_0)+I_1)\cap(l_1,f,g)$, and a general linear section of a reduced curve is reduced, we have that $I$ is also radical.

In order to show that the algebra defined by the ideal $I$ is level, let us start with the polynomial ring in nine variables 
\begin{equation*}
T=k[x,y,z_0,z_1,z_2,w_0,w_1,w_2,w_3],\end{equation*} and consider the ideal $J\subseteq T$ defined by 
\begin{equation*}
J = (xy,yf,yg,m_1,m_2,m_3),
\end{equation*}
where $m_1,m_2,m_3$ are the maximal minors of
\begin{equation*}\left(
\begin{array}{ccc}
  z_0&z_1&z_2\\
  w_0w_1&w_1w_2&w_2w_3\\
\end{array}
\right),\end{equation*}
and $f=m_1/w_2=z_1w_3-z_2w_1$, $g=m_3/w_1= z_0w_2-z_1w_0$.

A simple computation with \cite{Macaulay2} shows that the MFR of $J$ is:
\begin{equation*} 0 \rightarrow T^2(-6) \stackrel{\phi_2}{\rightarrow} T^6(-4)\oplus T(-5) \stackrel{\phi_1}{\rightarrow} T^5(-3)\oplus T(-2) \stackrel{\phi_0}{\rightarrow} J \rightarrow 0,\end{equation*}
where the maps $\phi_i$'s are given by the matrices: 
\begin{equation*}\phi_0=\left(\begin{array}{cccccc}
xy&yf&yg&m_1&m_2&m_3\end{array}\right),\end{equation*}
\begin{equation*}\phi_1=\left(
\begin{array}{ccccccc}
  0   & 0   & 0   & f  & g  & 0   & 0      \\
  w_2 & -w_0& 0   & -x & 0  & 0   & 0      \\
  0   &-w_3 & w_1 & 0  & -x & 0   & 0      \\
  -y  & 0   & 0   & 0  & 0  & z_0 & w_0w_1 \\
  0   & -y  & 0   & 0  & 0  & z_1 & w_1w_2 \\
  0   & 0   & -y  & 0  & 0  & z_2 & w_2w_3 \\
\end{array}
\right), 
\phi_2=
\left(
\begin{array}{cc}
  w_0w_1 & xz_0\\
  w_1w_2 & xz_1\\
  w_2w_3 & xz_2\\
  0      & g   \\
  0      & -f  \\
  0      & xy  \\
  y      & 0   \\
\end{array}
\right).
\end{equation*}

Thus, from the last module of the MFR, we deduce that the ideal $J$ is level of type two. Now we want to show that $I$ is obtained as a reduction of $J$ by a regular sequence.

Let $J^{'}$ be the extension of $J$ in the polynomial ring
\begin{equation*}R\otimes_kT=k[x_0,x_1,x_2,x_3,x,y,z_0,z_1,z_2,w_0,w_1,w_2,w_3].\end{equation*}\indent
Hence we have that $I$ is obtained from $J^{'}$ by reduction using the sequence
\begin{equation*}(x-l_0,y-l_1,z_0-C_0,z_1-C_1,z_2-C_2,w_0-Q_0,w_1-Q_1,w_2-Q_2,w_3-Q_3).\end{equation*}\indent
In order to prove that this sequence is regular for a general choice of the forms, it suffices to show that the algebra
$k[x_0,x_1,x_2,x_3]/I$ has dimension at most 1. Indeed, we have that $I$ contains the regular sequence
$(l_0l_1,Q_2f,Q_1g),$ as desired. Therefore, $R/I$ is the coordinate ring of a type two, level set of points of ${\mathbf P}^3$, since it is the reduction of a type two level algebra.

It now remains to check that the $h$-vector of $R/I$ is that of the statement. Notice that the algebra $A^{'}=R/(I,x_3)=k[x_0,x_1,x_2]/I^{'}$ is an artinian reduction of the algebra $R/I$ (since $x_3$ is a non-zero divisor of $R/I$), and therefore the $h$-vector $h=(h_0,h_1,...,h_e)$ and the graded Betti numbers $\beta_{i,j}$ are the same for the two algebras.

It is easy to check that the generators of $I$ (and hence those of $I^{'}$) have degrees
\begin{equation*}2,{\ }1+c_0+q_2= \left\lfloor\frac{e+3}{2}\right\rfloor,{\ }1+c_1+q_3=\left\lfloor\frac{e+4}{2}\right\rfloor ,{\ }c_0+q_1+q_2 = \left\lfloor\frac{2e+3}{3}\right\rfloor ,\end{equation*}\begin{equation*}c_0+q_2+q_3 =\left\lfloor\frac{2e+4}{3}\right\rfloor , {\ }c_1+q_2+q_3=\left\lfloor\frac{2e+5}{3}\right\rfloor .\end{equation*}\indent
The MFR and the $h$-vector of $A^{'}$ are related by the following well-known formula (e.g., see \cite[p. 131, point (j)]{FL}): 
\begin{equation}\label{sh}h(z)(1-z)^3=1+\sum_{i,j}(-1)^i\beta_{i,j}z^j,
\end{equation}
where we set $h(z)=\sum_{i=0}^eh_iz^i$.

From the maps $\phi_i$ described above, we have that the r.h.s. of (\ref{sh}) is\begin{equation*}
\begin{array}{rl}
1&-z^2-z^{1+c_0+q_2}-z^{1+c_1+q_3}-z^{c_0+q_1+q_2}-z^{c_0+q_2+q_3}-z^{c_1+q_2+q_3}\\ &+z^{1+c_0+q_1+q_2}+z^{1+c_0+q_2+q_3}
+z^{1+c_1+q_2+q_3}+z^{2+c_0+q_2}+z^{2+c_1+q_3}\\ &+z^{c_0+c_1+q_2+q_3}+z^{q_0+q_1+c_1+q_2+q_3}-2z^{e+3}\\ &=1-z^2-z^{\lfloor (e+3)/2\rfloor }-z^{\lfloor (e+4)/2\rfloor }-z^{\lfloor (2e+3)/3\rfloor }-z^{\lfloor (2e+4)/3\rfloor }\\ &-z^{\lfloor (2e+5)/3\rfloor }+z^{\lfloor (2e+6)/3\rfloor }+z^{\lfloor (2e+7)/3\rfloor }+z^{\lfloor (2e+8)/3\rfloor }\\ &+z^{\lfloor (e+5)/2\rfloor }+z^{\lfloor (e+6)/2\rfloor }+z^{e+1}+z^{e+2}-2z^{e+3}.
\end{array}
\end{equation*}

A standard (but tedious) computation shows that the entries of $h$-vector $h$ of the statement satisfy the equality (\ref{sh}), as we desired. This completes the proof of the theorem.
\end{proof}

\section{Non-unimodality for monomial level algebras}
We will next show the third result of this note: namely, as we said in the introduction, we will prove the existence of non-unimodal monomial (artinian) level algebras of codimension $r$, for any $r\geq 3$.

Remember that the {\it (degree) lexicographic order} is a total ordering, briefly indicated by $\lq \lq >"$, on the monomials of $S$ such that, if $P=y_1^{p_1}\cdot \cdot \cdot y_r^{p_r}$ and $Q=y_1^{q_1}\cdot \cdot \cdot y_r^{q_r}$ have the same degree, then $P>Q$ if and only if the first non-zero difference $p_i-q_i$ is positive.

\begin{ex}\label{ex4}
 Let us consider the inverse system module $M\subset S=k[y_1,y_2,y_3]$ generated by the last (according to the lexicographic order) 36 monomials of degree 12, namely $M=<y_1^2y_2^{10},y_1^2y_2^9y_3,...,y_3^{12}>$. An easy computation shows that the $h$-vector of $R/Ann(M)$ is \begin{equation*}(1,3,6,9,12,15,18,21,24,27,30,33,36).\end{equation*}\indent
Notice now that the form $F=y_1^6y_2^3y_3^3$ has all of its partial derivatives of order 1, 2 and 3 (which span vector spaces of dimension 3, 6 and 10, respectively) distinct from the partial derivatives of the forms generating $M$, since the latter forms and their derivatives are divisible by $y_1$ at most twice, whereas all the derivatives of order at most 3 of $F$ are divisible by $y_1$ at least thrice.
Therefore, the inverse system module $M^{'}=<M,F>$ generates an $h$-vector whose last four entries are: \begin{equation*}(27+10,30+6,33+3,36+1)=(37,36,36,37),\end{equation*} whence the $h$-vector $h$ of the codimension three artinian monomial level algebra $R/Ann(M^{'})$ is non-unimodal. In fact, we have \begin{equation*}h=(1,3,6,10,15,21,28,33,36,37,36,36,37).\end{equation*}
\end{ex}

\begin{rmk}\label{rmk}
i). The same construction can be used to produce infinitely many other examples. It suffices to generate $M$ by the last (according to the lexicographic order) $3e$ monomial of degree $e$, and to consider the form $F=y_1^{e-6}y_2^3y_3^3$. Then $M^{'}=<M,F>$ is, exactly as above, the inverse system module of a non-unimodal monomial level algebra of codimension 3 (its $h$-vector in fact ends with $(...,3e+1,3e,3e,3e+1)$).

ii). The above examples can be extended to construct non-unimodal monomial level algebras of any codimension $r\geq 3$. Indeed, if a monomial level algebra $k[x_1,x_2,x_3]/Ann(M)$ has $h$-vector $h=(1,3,h_2,...,h_e)$, then the monomial level algebra $k[x_1,...,x_r]/Ann(<M,y_4^e,...,y_r^e>)$ has $h$-vector $h^{'}=(1,r,h_2+r-3,....,h_e+r-3)$, which is clearly non-unimodal if $h$ is non-unimodal.
\end{rmk}

Hence we have shown:

\begin{thm}\label{thm5}
For every integer $r\geq 3$, there exist non-unimodal monomial artinian level algebras of codimension $r$.
\end{thm}

Theorem~\ref{thm5} immediately extends to ideals of points of the projective space, since monomial artinian algebras lift to reduced sets of points. Hence, we have supplied another proof of the following recent result of Migliore (\cite[Theorem 4.3]{Mi}):

\begin{cor}\label{cor6}
For every integer $r\geq 3$, there exist reduced level sets of points in ${\mathbf P}^r$ whose artinian reductions are non-unimodal.
\end{cor}

\section{Concluding questions}

We conclude our note with the following questions, that we regard as among the very next issues that now deserve to be addressed:

\begin{enumerate}
\item[i)] As far as the WLP is concerned, the existence problem for codimension 3 artinian level algebras, and for level sets of points of ${\mathbf P}^3$, remains now open only in the Gorenstein case: do all artinian Gorenstein algebras enjoy the WLP? (As we said in the introduction, \cite[Corollary 2.4]{HMNW}, provides a positive answer for the special case of complete intersections.) And if not, are there also reduced sets of points of ${\mathbf P}^3$ whose artinian reductions are Gorenstein and do not enjoy the WLP?

\item[ii)]As for monomial algebras, which is the largest type $t_0$ such that all the codimension 3 monomial level algebras of type $t\leq t_0$ are unimodal? It is easy to see that $t_0>0$, whereas here we have shown that $t_0$ is finite (in particular, $t_0\leq 36$).

\item[iii)] Similarly to what the second author has done in \cite[Remark 5, iii)]{Za} for arbitrary level algebras, can we extend Example~\ref{ex4} to $h$-vectors with three maxima? With any number of maxima?
\end{enumerate}

\bibliographystyle{amsalpha} 

\end{document}